\documentclass[a4paper,11pt]{amsart}
\usepackage[latin1]{inputenc}
\usepackage{amsmath,amsfonts,amssymb}
\def\R{\mathbb{R}}
\def\C{\mathbb{C}}

\def\B{\mathbb{B}}
\def\H{\mathbb{H}}

\def\Re{\mathrm{Re}\,}

\def\tsp{{\,}^t\!}
\def\vp{\varphi}
\def\ve{\varepsilon}
\def\D{\Delta}

\def\dist{\mathrm{dist}\,}
\def\pri{\,'\!}

\theoremstyle{plain}
\newtheorem{theo}{Theorem}
\newtheorem{cor}{Corollary}
\newtheorem{prop}{Proposition}
\newtheorem{lm}{Lemma}
\newtheorem{defi}{Definition}
\newenvironment{dem}{\noindent{\bf{Proof}} \\}{\hfill $\square$ \newline}
\newenvironment{pr}{\begin{trivlist}\item[]\textit{Proof.}}{\hfill $\square$\end{trivlist}}

\theoremstyle{remark}
\newtheorem{rem}{Remark}

\title[Proper pseudo-holomorphic maps]{Proper pseudo-holomorphic maps between strictly pseudoconvex regions}
\author{Léa Blanc-Centi}
\address{Universités de Marseille, Université de Provence, L.A.T.P.\\ 39 rue Joliot-Curie\\ 13453 Marseille Cedex 13\\ FRANCE}
\email{lea@cmi.univ-mrs.fr}

\begin{document}

\begin{abstract}
We prove that proper pseudo-holomorphic maps between strictly pseudoconvex regions in almost complex manifolds extend to the boundary. The key point is that 
the Jacobian of such a map is far from zero near the boundary, and the proof is mainly based on an almost complex analogue of the scaling method. 
We also establish the link between the regularity of the extension and the regularity of the almost complex structures, and we determine explicit estimates for the Hölderian norms. As a corollary, we get in the smooth case a necessary and sufficient condition on the almost complex structure of the target's space for the smooth extension of proper pseudo-holomorphic maps. 
\end{abstract}

\maketitle

\section{Introduction}
The regularity up to the boundary of a proper holomorphic map between strictly pseudoconvex bounded domains $D$ and $D'$ of $\C^n$ has been widely studied when $n\ge 2$, and is now as well understood as the one-dimensional case.  
A continuous map $F:D\to D'$ is said to be {\em proper} if $F^{-1}(K)$ is compact for every compact set $K$ in $D'$.
If $D$ and $D'$ have $\mathcal{C}^r$-boundaries ($r\ge 2$), such a map $F$ has a $\mathcal{C}^{r-1/2}$-extension to the boundary, and this is the maximal regularity \cite{Pinchuk,Co} that can be expected. Various authors have contributed to this result. 
We just mention Ch. Fefferman, who proved in 1974 that if $D$ and $D'$ have smooth boundaries, and if $F$ is a biholomorphism, then it extends smoothly to the boundary \cite{Feff}. We refer to the survey of F. Forstneri\v{c} \cite{Forst} for a thorough history. 
 
Since an almost complex structure is generically non-integrable, a natural question is to ask whether this result is always valid in the almost complex situation. In order to establish an analogue of the previous result, we will need to adapt objects and tools specific to the integrable case. Note for example that there is no more notion of an analytic set, and that the Jacobian of a pseudo-holomorphic map is not pseudo-holomorphic. 
There are various approaches for proving theorems of extension. Originally introduced by E. Bishop, the use of analytic discs in order to relate the behavior at the boundary and the behavior in the domain presents the advantage of being quite easily transposed to the almost complex case.  
The method based on analytic discs has provided geometric proofs of various versions of Fefferman's theorem \cite{L1,Tumanov}, even in the almost complex case \cite{CGS,GSmodele}.\\

Our aim here is to study the boundary behavior of proper pseudo-holomorphic maps between strictly pseudo-convex regions. More precisely, we consider the following situation.
Let $D$ be a bounded domain in some smooth (real) manifold, and $J$ be an almost complex structure of class $\mathcal{C}^1$ on $\overline{D}$, smooth in $D$. Throughout this paper, we will say that $(D,J)$ is a {\em strictly pseudoconvex region} if $D$ is defined by $\{\rho<0\}$, with $\rho$ being a $\mathcal{C}^2$-regular defining function of $\partial D$, strictly $J$-plurisubharmonic on $\overline{D}$. We will call $(D,J)$ a {\em strictly pseudoconvex region of class $\mathcal{C}^r$} when $\rho$ and $J$ are at least of class $\mathcal{C}^r$.
In the complex situation, the regularity is thus the regularity of the boundary.

The first step is to get the Hölder $1/2$-continuous extension, which comes from a sized estimate of the set of regular values, and from estimates of the Kobayashi metric. To obtain more regularity, the main obstacle compared with the biholomorphic case is obviously the existence of critical points. Thus we begin with studying the locus of all these points.
We recall that in the complex case, S. Pinchuk has showed using the scaling method that every proper holomorphic map between two bounded strictly pseudoconvex domains in $\C^n$, with $\mathcal{C}^2$-boundaries, is locally biholomorphic (\cite{P1}, see also \cite{Pinchuk}). 
 We prove that this statement is always true in the almost complex case, at least near the boundary:

\begin{theo}\label{theo:jacobiencompact} 
Let $(D,J)$ and $(D',J')$ be some strictly pseudoconvex regions of the same dimension, and $F$ be a proper pseudo-holomorphic map from $D$ to $D'$. Then the Jacobian of $F$ is far from zero near the boundary of $D$.
\end{theo}

In particular, the critical points of $F$ form a compact set in $D$. Our proof is based on an almost complex adaptation of the scaling method (see also \cite{CGS}). The idea of this method consists in dilating anisotropically the domains, in order to construct a limit map between simple model domains. When the manifolds are almost complex, the transformations operating on the domains are not pseudo-holomorphic, so we simultaneously dilate the almost complex structures. Note that, if $n=2$, one can normalize the initial structures to obtain the standard structure as a limit. In the general case, the limit almost complex structures are not necessarily integrable.
 
The essential tool to get more regularity for the extension to the boundary is the study of some family of pseudo-holomorphic discs attached to the boundary of the domain, which is invariant under the action of pseudo-biholomorphisms. We obtain that the regularity of the extension depends on the regularity of the almost complex structures at the boundary, and we also give explicit estimates of the Hölderian norms, by using the results of \cite{3}. More precisely:

\begin{theo}\label{theo:estimationpropre}
Let $(D,J)$ and $(D',J')$ be strictly pseudoconvex regions of the same dimension, respectively $\mathcal{C}^r$ and $\mathcal{C}^{r'}$ where $r,r'\ge 2$ are not integers.
Then every proper pseudo-holomorphic map $F:D\to D'$ has a $\mathcal{C}^{s}$-extension to $\overline{D}$, where 
\begin{itemize}
\item $s=\mathrm{min}\,(r,r')$ if $|r'-r|\ge 1$;
\item $s=\mathrm{max}\,(r-1,r'-1)$ if $|r'-r|<1$. 
\end{itemize}
Moreover, 
$$||F||_{\mathcal{C}^{s-1}(\bar{D})}\le c(s)||(F,\tsp(dF)^{-1})||_\infty\left(1+\frac{c'}{\sqrt{\lambda_{N^*M'}^{J'}}}\right).$$
\end{theo}

Here, $\lambda_{N^*M'}^{J'}$ denotes the smallest eigenvalue of the Levi form of some function defined by means of the equations of the conormal bundle $N^*M'$. We will call it the {\em minimal curvature of $N^*M'$}. 

When $(D,J)$ and $(D',J')$ are smooth, we thus obtain that $F$ has a smooth extension. 
Note that $J'=F_*J$ near the boundary by Theorem \ref{theo:jacobiencompact}. Hence, we get, as in the biholomorphic case \cite{CGS} a necessary and sufficient condition on $J'$ for the smooth extension of $F$:

\begin{cor}\label{cor:CNS}
Let $(D,J)$ and $(D',J')$ be strictly pseudoconvex regions of the same dimension, with $(D,J)$ of class $\mathcal{C}^\infty$. A proper pseudo-holomorphic map $F:D\to D'$ extends smoothly to the boundary if and only if $\partial D'$ is smooth and $J'$ extends smoothly on $\overline{D'}$.
\end{cor}

This paper is organized as follows. The second section is made of some preliminaries about almost complex manifolds. In the third section, we prove the boundary distance preserving property for a proper pseudo-holomorphic map, which leads to the existence of the Hölder $1/2$-continuous extension. Section 4 is devoted to the application of the scaling method; we then study the properties of the limit map. Finally, we prove Theorems \ref{theo:jacobiencompact} and \ref{theo:estimationpropre} and Corollary \ref{cor:CNS} in the last section. 



\section{Preliminaries}
\subsection{Strictly pseudoconvex domains in almost complex manifolds}
Let us recall some definitions.

\begin{defi}
An {\em almost complex structure} on a smooth (real) manifold $M^{2n}$ is a $(1,1)$-tensor $J$, that is, a section from $M$ to $End(TM)$, such that $J^2=-Id$. .
\end{defi}

Every almost complex structure admits a Hermitian metric, and provides an orientation on the manifold.

\begin{defi}
A map $F:(M,J)\to(M',J')$ of class $\mathcal{C}^1$ between two almost complex manifolds is said to be {\em $(J,J')$-holomorphic} if $J'\circ dF=dF\circ J$.
\end{defi}

If $(M,J)$ is the unit disc of $\C$ (that is, $\D\subset\R^2$ equipped with the standard complex structure), we say that $F$ is a {\em $J$-holomorphic disc}.
A. Nijenhuis and W. Woolf \cite{NW} have proved that there exist such maps. Moreover, one can prescribe $F(0)$ and $dF_0(\partial/\partial\, x)$.

As in the complex situation, the maps verifying the equation of pseudo-holomorphy inherit their smoothness from the smoothness of the almost complex structures: if $J$ and $J'$ are of class $\mathcal{C}^r$, then every $(J,J')$-holomorphic map  is of class $\mathcal{C}^{r+1}$.
We also have:

\begin{lm}\label{lm:ConservationOrientation}
Let $F$ be a pseudo-holomorphic map between almost complex manifolds of the same dimension. Then either $F$ preserves or $F$ inverses the orientation provided by the almost complex structures.
\end{lm}

\begin{pr}
Once having fixed local coordinates, we have to show that the sign of the Jacobian of $F$ is constant on $D$.
The almost complex structure $J$ verifies $J_p=PJ_{st}P^{-1}$ for some matrix $P$, where $J_{st}=\left(\begin{array}{cc}0&-I_n\\I_n&0\end{array}\right)$ is the standard complex structure on $\R^{2n}$. Such a factorization is not unique, but if $J_p=PJ_{st}P^{-1}=QJ_{st}Q^{-1}$, then the matrix $Q^{-1}P$ commutes with $J_{st}$ and hence $\mathrm{det}(Q^{-1}P)>0$. Thus the sign of $\mathrm{det}P$ only depends on $p$ and $J$. We denote it by $\delta^J(p)$. 

Given $p_0\in D$, let $(e_1^{(p_0)},\hdots,e_{2n}^{(p_0)})$ be a basis such that $\forall k\ge 1,\ e_{n+k}^{(p_0)}=J_{p_0}e_k^{(p_0)}$. Denote by $P_p$ the matrix of $(e_1^{(p_0)},\hdots,e_{n}^{(p_0)},J_pe_1^{(p_0)},\hdots,J_pe_{n}^{(p_0)})$, and let $V$ be a neighborhood of $p_0$ in which the matrix $P_p$ remains invertible. In particular, for every $p\in V$, $\mathrm{det}P_p$ and $\mathrm{det}P_{p_0}$ have the same sign.
It follows that $\delta^J$ is locally constant on $D$, hence constant.

For any $p\in D$ and $q\in D'$, we write $J_p=P_pJ_{st}P_p^{-1}$ and $J'_q=P'_qJ_{st}{P'}_q^{-1}$. The $(J,J')$-holomorphy of $F$ gives:
$\forall p\in D,\ dF_pJ_p=J'_{F(p)}dF_p$, so $({P'}_{F(p)}^{-1}dF_pP_p)$ commutes with the complex standard structure and $\mathrm{det}({P'}_{F(p)}^{-1}dF_pP_p)>0$. Finally, the sign of the Jacobian of $F$ is equal to $\delta^J\times\delta^{J'}$ at any non-critical point.
\end{pr}

\begin{defi}\label{def:strpsh}
Let $(M,J)$ be an almost complex manifold and $\rho$ be a $\mathcal{C}^2$-smooth function from $M$ to $\R$. For all $X\in TM$, define $d^c_J\rho(X)=-d\rho(JX)$ and $\mathcal{L}^J\rho=d(d^c_J\rho)(X,JX)$. The quadratic form $\mathcal{L}^J\rho$ is called the {\em Levi form} of $\rho$.
The function $\rho$ is said to be {\em strictly $J$-plurisubharmonic} if its Levi form is positive definite.
\end{defi}

\begin{rem}\label{rem:CNSstrpshC2}
One may show that a $\mathcal{C}^2$-regular map $\rho:\Omega\to\R$ is strictly $J$-plurisubharmonic if and only if all the compositions $\rho\circ h$, for any $J$-holomorphic disc $h$ valued in $\Omega$, are strictly subharmonic on the unit disc $\Delta$. Indeed, there is the following link between the Laplacian of $\rho\circ h$ and the Levi form of $\rho$ \cite{ref7Lee, Haggui}: 
$$\forall\zeta\in\D,\ \bigtriangleup(\rho\circ h)_{\zeta}=\mathcal{L}^J_{h(\zeta)}(\rho)\left(\frac{\partial h}{\partial x}(\zeta)\right).$$
\end{rem}

\subsection{Holomorphic mappings between simple model domains}
Model almost complex structures naturally appear as limits of rescaled almost complex structures. We refer to \cite{GSmodele} for a detailed treatment of model structures.

Throughout this paper, we denote by $(x_0,y_0,\hdots,x_n,y_n)$ the coordinates in $\R^{2n+2}$, and by $z=(z_0,\hdots,z_n)=(z_0,\pri z)$ the associated complex coordinates. So $\R^{2n+2}$ may be identified with $\C^{n+1}$ by means of the standard complex structure
$$\mathcal{J}^{(n+1)}_{st}=\left(\begin{array}{ccccc}0&-1& & & \\ 1&0& &(0) & \\ & &\ddots& & \\ &(0) & &0&-1\\ & & &1&0\end{array}\right).$$

\begin{defi}\label{coeffB}
An almost complex structure $J$ on $\R^{2n+2}$ is called a {\em model structure} if it is defined by
\begin{equation}\label{formemodele}
J(z)=\left(\begin{array}{cc}\mathcal{J}_{st}^{(1)}&B^J(\pri z)\\ 0&\mathcal{J}_{st}^{(n)}\end{array}\right)\,,
\end{equation}
where $B^J(\pri z)\in\mathcal{M}_{2,2n}(\R)$ is $\R$-linear in $x_1,\hdots,x_n,y_1,\hdots,y_n$. A pair $(\Sigma,J)$ is called a {\em model domain} if 
$\Sigma=\{z\in\C^{n+1}/\ \Re z_0+P(\pri z,\pri \bar{z})<0\}$, where $P$ is some real homogeneous polynomial of degree 2, and $J$ is a model structure such that $\Sigma$ is strictly $J$-pseudoconvex at 0.  
\end{defi}

\noindent For such a matrix $B^J$, the corresponding complex matrix is
$$B_{\C}('z)=\left(B_{2k-1}^{2j-1}+iB_{2k-1}^{2j}\right)_{1\le j,k,\le n}=\left(\sum_{k=1}^n(a_{1,k}^Jz_k+b_{1,k}^J\bar{z}_k)\ \hdots\ \sum_{k=1}^n(a_{n,k}^Jz_k+b_{n,k}^J\bar{z}_k)\right)$$
where $a_{j,k}^J$ and $b_{j,k}^J$ are complex constants.

\begin{defi}
The model structure $J$ given by (\ref{formemodele}) is {\em simple} if for all $j,k$, $a_{j,k}^J=0$.
\end{defi}

For any model domain $(\Sigma,J)$, there exists a simple model structure $\mathcal{J}$ and a $(J,\mathcal{J})$-biholomorphism between $\Sigma$ and $\H$ fixing $(-1,\pri 0)$, where 
$\H=\{z\in\C^{n+1}/\ \Re z_0+||\pri z||^2<0\}$
is the Siegel half-plane \cite{Leebis}.

The proof shows that the map constructed is in fact a global diffeomorphism of $\C^{n+1}$, given by
$$\Psi(z)=(z_0+\psi(\pri z),\pri\Psi(\pri z))$$
where $\psi$ is a polynomial of degree at most two in $\pri z$ and $\pri\bar{z}$, and $\pri\Psi:\C^n\to\C^n$ is $\C$-linear in $\pri z$. Consequently:
\begin{enumerate}
\item for all $t<0$, $\Psi(t,\pri 0)=(t,\pri 0)$;
\item if the model domain $\Sigma$ is given by the equation $0=\tilde{r}(z)=\Re z_0+P(\pri z,\pri \bar{z})<0$ with $P$ some real homogeneous polynomial of degree two on $\C^n$, then 
$$\tilde{r}\circ\Psi^{-1}(z)=\Re z_0+||\pri z||^2.$$
\end{enumerate}
The proof also implies that if the model structure $\mathcal{J}$ is integrable, we can in fact prescribe $\mathcal{J}=\mathcal{J}_{st}$.\\

Because of the special form of simple model structures, a pseudo-holomorphic map  $F=(F_0,\pri F)$ between simple model domains has an interesting behavior. We can suppose that the domains are both the Siegel half-plane $\H$. 

\begin{prop}{\bf \cite{Leebis}}
Assume that $\mathcal{J}$ and $\mathcal{J}'$ are \textup{non-integrable} simple model structures on $\H$. If $F:(\H,\mathcal{J})\to(\H,\mathcal{J}')$ is a pseudo-holomorphic map, then there exists a real constant $c$ such that
\begin{equation}\label{FormeCasNonInt}
\forall z=(z_0,\pri z)\in\H,\ F(z)=(cz_0+f(\pri z),\pri F(\pri z)),
\end{equation}
where $f:\C^n\to\C$ is anti-holomorphic and $\pri F:\C^n\to\C^n$ is holomorphic (with respect to the standard complex structure).
\end{prop}

\begin{pr}
The proof is given in \cite{Leebis} in the case of a pseudo-biholomorphism. We sketch it in the general case, using the notations of Definition \ref{coeffB}.

The computation of the coefficients of the Nijenhuis tensor $N_\mathcal{J}$ and the hypothesis that $\mathcal{J}$ is non-integrable give some $j,k$ such that $b^\mathcal{J}_{j,k}-b^\mathcal{J}_{k,j}\not=0$. Let us fix such a pair $(j,k)$. We may identify the coefficients in $dF\left(N_\mathcal{J}\left(\frac{\partial}{\partial{z}_j},\frac{\partial}{\partial{z}_k}\right)\right)=N_{\mathcal{J}'}\left(dF\left(\frac{\partial}{\partial{z}_j}\right),dF\left(\frac{\partial}{\partial{z}_k}\right)\right)$, and get equivalently
\begin{eqnarray}\label{F0Z0}
 & &\forall l\ge 1,\ \frac{\partial F_l}{\partial z_0}=\frac{\partial F_l}{\partial \bar{z}_0}=0\label{indepz0}\\
 & &\frac{\partial F_0}{\partial z_0}=\frac{1}{b_{j,k}^\mathcal{J}-b_{k,j}^\mathcal{J}}\,\sum_{l,m=1}^n\frac{\partial\bar{F}_l}{\partial\bar{z}_j}\frac{\partial\bar{F}_m}{\partial\bar{z}_k}(b^{\mathcal{J}'}_{l,m}-b^{\mathcal{J}'}_{m,l})\,.
\end{eqnarray}
The system obtained in (\ref{indepz0}) means that $F_1,\hdots, F_n$ do not depend on $z_0$ and $\bar{z}_0$ (precisely, of $x_0$ and $y_0$). Moreover $F$ is $(\mathcal{J},\mathcal{J}')$-holomorphic, and since $\mathcal{J},\mathcal{J}'$ are simple model structures, we have
$$dF\circ\mathcal{J}=\mathcal{J}'\circ dF\Longrightarrow\left\{\begin{array}{l}
\mathcal{J}_{st}^{(1)}\circ d_{z_0}F=d_{z_0}F\circ\mathcal{J}_{st}^{(n)}\\
\mathcal{J}_{st}^{(n)}\circ d(\pri F)=d(\pri F)\circ\mathcal{J}_{st}^{(n)}
\end{array}\right.\,.$$
So $\pri F:\pri z\mapsto (F_1(\pri z),\hdots, F_n(\pri z))$ and $z_0\mapsto F_0(z_0,\pri z)$ are holomorphic (with respect to the standard complex structure). This implies that $\frac{\partial F_0}{\partial z_0}$ is independent of $z_0$, and anti-holomorphic in $\pri z$: 
$$F_0(z)=c(\pri z)z_0+f(\pri z),$$
where $c$ and $f$ are anti-holomorphic.
For any $\pri z\in\C^n$, the map $z_0\mapsto F_0(z_0,\pri z)$ is defined on $\{\zeta\in\C/\ \Re\zeta<-||\pri z||^2\}$ and takes its values in $\{\zeta\in\C/\ \Re\zeta<-||\pri F(\pri z)||^2\}$. This needs $c$ to be real valued and anti-holomorphic, hence constant.
\end{pr}

\section{Hölder $1/2$-continuous extension}
\subsection{Regular values of proper maps}
We recall that a \textit{critical point} of a $\mathcal{C}^1$ map $F$ from $D$ to $D'$ is a point $p\in D$ such that the Jacobian of $F$ vanishes at $p$. In the sequel, we will denote by $C$ the set of all critical points. A {\em critical value} is the image by $F$ of some critical point. Every point of $D'$ that is not a critical value, even if it is not in $F(D)$, is called a {\em regular value}. 

\begin{rem}\label{dimH}
If $F:(D,J)\to(D',J')$ is pseudo-holomorphic, then for any critical point $p\in D$, the subspace $\mathrm{Ker}\,dF_p$ contains a subspace of dimension 2 since it is preserved by $J_p$. According to Theorem 3.4.3 in \cite{Fe}, we get that the Hausdorff dimension of the set $F(C)$ of all critical values is less than $2n-2$.
\end{rem}

We also show the following property of proper pseudo-holomorphic maps:

\begin{lm}\label{lm:valeursregulieres}
Let $(D,J)$ and $(D',J')$ be strictly pseudoconvex regions of dimension $2n$, and $F$ be a proper pseudo-holomorphic map from $D$ to $D'$. Then there is at least one point $p\in D$ such that $F(p)$ is a regular value.
\end{lm} 

\begin{pr}
It suffices to prove that the open set $D\setminus C$ is not empty. Indeed, this will imply, by means of the rank theorem, that the (Hausdorff) dimension of $F(D\setminus C)$ is equal to $2n$. According to Remark \ref{dimH}, we will get $F(C)\not= F(D\setminus C)$.\\
Assume by contradiction that $C=D$. Then the maximal rank $r_0$ of $dF$ in $D$ is less than $2n-2$, and non-zero since $F$ is necessarily non-constant. The locus $C_0$ where the rank of $dF$ is $r_0$ is an open set, therefore $F(C_0)$ is a submanifold of (Hausdorff) dimension $r_0$.

As the rank of the map $dF$ is less than $r_0-1$ (in fact, $r_0-2$) on $D\setminus C_0$, Theorem 3.4.3 in \cite{Fe} gives $\mathrm{dim}_HF(D\setminus C_0)\le r_0-1$. Thus $F(C_0)\not= F(D\setminus C_0)$.
This allows us to pick $q\in F(C_0)\setminus F(D\setminus C_0)$. The set $N=F^{-1}(\{q\})$ included in the open set $C_0$ is then a submanifold of $D$ of dimension $2n-r_0$. Moreover, $N$ is compact, and its tangent bundle $\mathrm{Ker}\,dF$ is preserved by $J$. Hence $J$ induces an almost complex structure on $N$.

Let $\rho$ be some $\mathcal{C}^2$-regular and strictly J-plurisubharmonic function on $D$.
The map $\rho$ reaches its maximum on $N$ at some point $p$. According to \cite{NW}, there exists a pseudo-holomorphic disc $h$ which takes its values in $N$, centered at $p$ and verifying $\frac{\partial h}{\partial x}(0)\not=0$. By the maximum principle, the strictly subharmonic function $\rho\circ h$ is constant. This leads to a contradiction since $0=\bigtriangleup(\rho\circ h)_{p}=\mathcal{L}^J_{p}(\rho)\left(\frac{\partial h}{\partial x}(0)\right)>0$ in view of Remark \ref{rem:CNSstrpshC2}.
\end{pr}

The next statement is well-known in the complex case:

\begin{prop}\label{prop:surjectivite}
Let $(D,J)$ and $(D',J')$ be strictly pseudoconvex regions of the same dimension. If $F$ is proper and pseudo-holomorphic from $D$ to $D'$, then $F(D)=D'$. Moreover, all the regular values of $F$ have the same (finite) number of antecedents, and they form a connected-by-path open set that is dense in $D'$.
\end{prop}

\begin{pr}
Let us recall (see for example \cite{Spivak}) that if $X$ and $Y$ are two oriented connected manifolds of the same dimension, and $F:X\to Y$ is a smooth proper map, then for any regular value $q\in Y$, the degree of $F$ is equal to $\sum_{p\in F^{-1}\{q\}}\mathrm{sgn}(\mathrm{det}\,dF_p)$. In particular, if $q$ is not in the image of $F$, then the degree of $F$ is
zero.
Let $q\in F(D)\setminus F(C)$ be a regular value of $F$. The degree of $F$ is either positive or negative by Lemma \ref{lm:ConservationOrientation}, thus $F(D)=D'$ and the preimage of any regular value of $F$ has exactly $\mathrm{deg}\,F$ elements.

The set of all critical points of $F$ is closed in $D$. Indeed, $F$ is proper hence closed, and $D'\setminus F(C)$ is a open set in $D'$. According to Proposition 14.4.2. in \cite{Rudin} and Proposition \ref{dimH}, $F(C)$ has no interior. The fact that $D'\setminus F(C)$ is connected-by-path comes from classical geometric arguments, as in the complex case.
\end{pr}

\subsection{Boundary distance preserving property}
The aim of this subsection is to prove the following statement, which will give as a corollary the Hölder $1/2$-continuous extension of a proper pseudo-holomorphic map:

\begin{prop}\label{conservationdist}
Let $(D,J)$ and $(D',J')$ be strictly pseudoconvex regions of dimension $2n$, and $F$ be a proper pseudo-holomorphic map from $D$ to $D'$. Then there exist $c_1,\ c_2>0$ such that
$$\forall p\in D,\ c_1\le\frac{\mathrm{dist}(F(p),\partial D')}{\mathrm{dist}(p,\partial D)}\le c_2.$$
\end{prop}

\begin{pr}
The Hopf lemma in the almost complex situation \cite{CGS} applies to the map $\rho'\circ F$ on $D$ and gives some $c'>0$ such that$$\forall p\in D,\ |\rho'(F(p))|\ge c'\,\mathrm{dist}(p,\partial D).$$
The idea of the proof of Proposition \ref{conservationdist} is to construct a $J$-plurisubharmonic map on $D'$, say $\rho'\circ F^{-1}$, in order to apply again the Hopf lemma. More precisely, we define 
$$u:D'\ni q\mapsto\underset{p\in F^{-1}(\{q\})}{\mathrm{Max}}\{\rho(p)\}.$$ 
Note that $u$ takes negative values, and is continuous on $D'\setminus F(C)$. Indeed, for any $q$ in $D'\setminus F(C)$, the compact set $K=F^{-1}(\{q\})$ consists in non critical, hence isolated, points, so is finite. Set $K=\{p_1,\hdots,p_k\}$, where $k=|\mathrm{deg}\,F|$. According to the inverse function theorem, 
  one can construct for each $j=1,\hdots,k$ a neighborhood $V_j$ of $p_j$, and a neighborhood $W$ of $q$ included in the open set $ D'\setminus F(C)$ such that $F$ induces a $\mathcal{C}^1$-diffeomorphism from $V_j$ to $W$, and 
$F^{-1}(W)=\bigsqcup V_j$. 
For $j=1,\hdots,k$, set $F_j=F_{|V_j}:V_j\to W$ and $u_j=\rho\circ F_j^{-1}$. The maps $u_j$ are continuous and $J$-plurisubharmonic. Consequently, $u=\underset{1\le j\le k}{\mathrm{Max}}\,u_j$ is continuous. 

We also get that $u$ is $J$-plurisubharmonic on $W'$, in the following sense: the composition $u\circ h$ is subharmonic on the unit disc $\Delta\subset\C$ for every $J$-holomorphic disc $h:\Delta\to W$ (this extends Definition \ref{def:strpsh} to the case of upper semicontinuous function). Thus $u$ is locally $J$-plurisubharmonic on $D'\setminus F(C)$, which is equivalent to be globally $J$-plurisubharmonic.

Since $\mathrm{dim}_H(F(C))\le 2n-2$, we obtain that $\mathrm{limsup}(u)$ is plurisubharmonic on the whole $D'$ (see \cite{CS}). By the Hopf lemma, there is some constant $c>0$ such that for any $q\in D'$, $|\mathrm{limsup}(u)(q)|\ge c\,\mathrm{dist}(q,\partial D')$. The map $\rho$ being continuous,
$$\forall p\in D,\ |\rho(p)|\ge |\mathrm{limsup}(u)(F(p))|\ge c\,\mathrm{dist}(F(p),\partial D'),$$
and we obtain the desired inequalities.
\end{pr}

\begin{cor}\label{cor:prolC1/2}
Let $(D,J)$ and $(D',J')$ be strictly pseudoconvex regions of dimension $2n$, and $F$ be a proper pseudo-holomorphic map from $D$ to $D'$. Then $F$ has a continuous extension $\bar{F}:\overline{D}\to\overline{D'}$ such that $\bar{F}(\partial D)\subset\partial D'$. Moreover, $\bar{F}$ is H\"older continuous with exponent 1/2.
\end{cor}

\begin{pr}
Once $F$ is known to preserve the distance to the boundary, Corollary \ref{cor:prolC1/2} comes from estimates of the infinitesimal Kobayashi pseudometric, exactly as in the proof of Proposition 3.3 in \cite{CGS}. We sketch the arguments here, for the sake of completeness. For every $p\in D$ and $v$ a tangent vector at point $p$, we set
$$K_{(D,J)}(p,v)=\text{inf}\{\alpha>0/\ \exists h\in\mathcal{O}^J(\Delta,D)\ \text{with}\ h(0)=p\ \text{and}\ (\partial h/\partial x)(0)=v/\alpha\},$$
which is well-defined according to \cite{NW}. 
Let us recall the following result \cite{GS,CGS}: under our hypotheses, there exists a constant $C>0$ such that
$$\forall p\in D,\ \forall v\in T_pM,\ (1/C)\frac{||v||}{\mathrm{dist}(p,\partial D)^{1/2}}\le K_{(D,J)(p,v)}\le C\frac{||v||}{\mathrm{dist}(p,\partial D)}.$$
 Then, because of the decreasing property of the infinitesimal Kobayashi pseudometric, we have for any $p\in D$ and $v\in T_pM$:
$$C_1\frac{||dF_p(v)||}{\mathrm{dist}(F(p),\partial D')^{1/2}}\le K_{(D',J')}(F(p),dF_p(v))\le K_{(D,J)}(p,v)\le C_2\frac{||v||}{\mathrm{dist}(p,\partial D)}$$
which implies the estimate
$\displaystyle |||dF_p|||\le C\frac{||v||}{\mathrm{dist}(p,\partial D)^{1/2}}$.
This gives the statement by Hardy-Littlewood's theorem.
\end{pr}

\section{The scaling method}
The proof of Theorem \ref{theo:jacobiencompact} is based on the scaling method. To make the transformations holomorphic, we also rescale the almost complex structures as in \cite{GSmodele} (see also \cite{Leebis}).

The idea is to rectify $\partial D$ and $\partial D'$ by means of successive changes of variable, in order to obtain $\partial\H$ (to make the transformations holomorphic, we also rescale the almost complex structures):

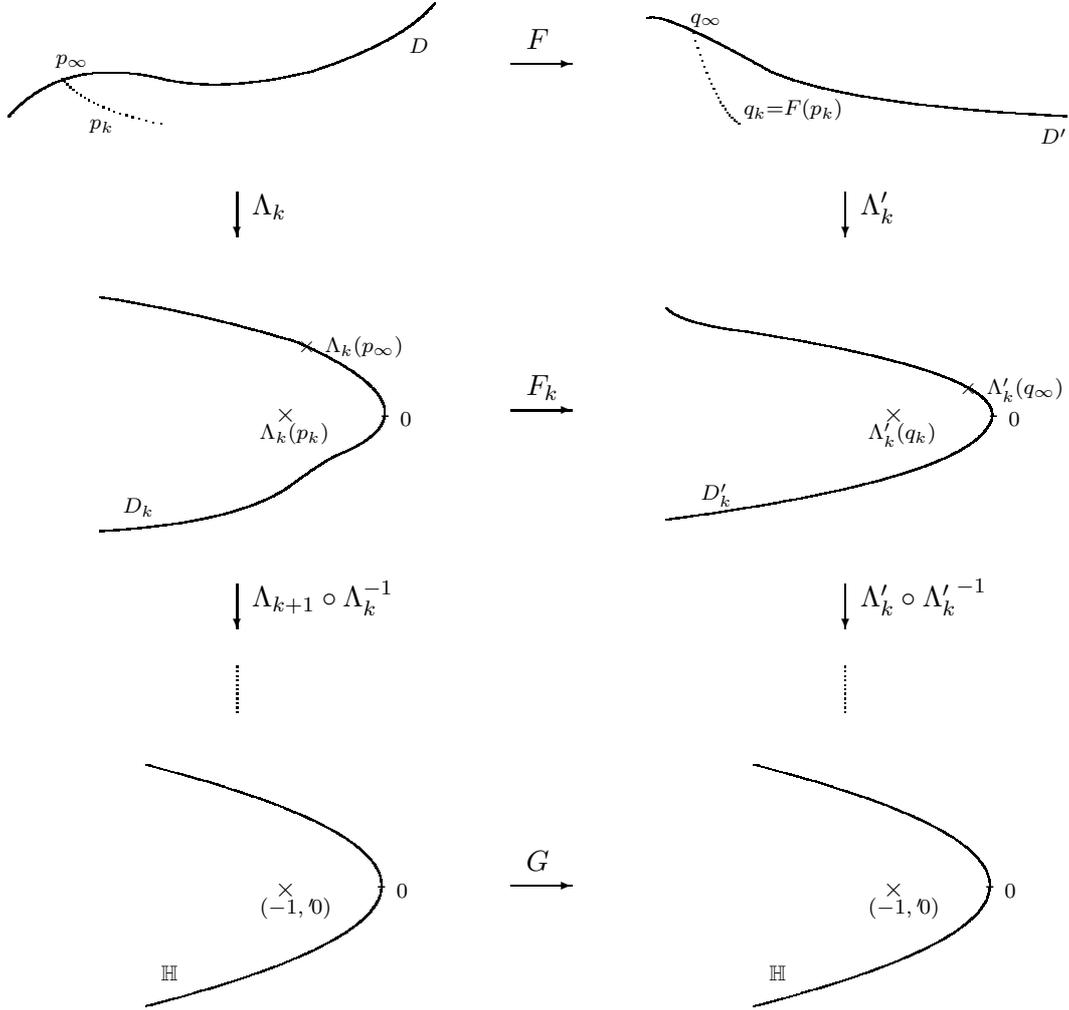
\begin{picture}(14,14)
\qbezier(0,11.8)(0.7,12.6)(2,12.3)
\qbezier(2,12.3)(2.8,12.1)(4,12.4)
\qbezier(4,12.4)(5.2,12.8)(5.6,13.3)
\qbezier[20](0.7,12.3)(1,11.9)(2,11.7)
\put(0.55,12.55){$\,_{p_\infty}$}
\put(1,11.7){$\,_{p_k}$}
\put(5.2,12.7){$\,_{D}$}

\put(6.6,12.5){\vector(1,0){0.8}}
\put(6.8,12.7){$F$}

\qbezier(8.4,13.1)(8.6,13.2)(10,12.4)
\qbezier(10,12.4)(11,11.95)(13.9,11.8)
\qbezier[20](9,12.9)(9.3,11.9)(9.6,11.7)
\put(8.9,13.1){$\,_{q_\infty}$}
\put(9.6,11.9){$\,_{q_k=F(p_k)}$}
\put(13.5,11.5){$\,_{D'}$}

\put(3,10.8){\vector(0,-1){0.6}}
\put(3.2,10.5){$\Lambda_k$}
\put(11,10.8){\vector(0,-1){0.6}}
\put(11.2,10.5){$\Lambda'_k$}

\qbezier(1.2,9.4)(2.5,9.2)(3.8,8.8)
\qbezier(3.8,8.8)(5.8,7.9)(4.3,7.3)
\qbezier(3.7,6.9)(4.1,7.2)(4.3,7.3)
\qbezier(1.2,6.3)(3,6.4)(3.7,6.9)
\put(3.5,7.75){{\small $\times$}}
\put(3.3,7.6){$_{\Lambda_k(p_k)}$}
\put(3.8,8.75){$_{\times\ \Lambda_k(p_\infty)}$}
\put(1.5,6.6){$_{D_k}$}
\put(4.9,7.75){{\small -} $_0$}

\put(6.6,7.9){\vector(1,0){0.8}}
\put(6.8,8.1){$F_k$}

\qbezier(8.65,9.25)(8.85,9.05)(9.65,8.95)
\qbezier(9.65,8.95)(13.45,8.35)(12.85,7.65)
\qbezier(8.65,6.45)(12.35,6.95)(12.85,7.65)
\put(11.5,7.75){{\small $\times$}}
\put(11.3,7.6){$_{\Lambda'_k(q_k)}$}
\put(12.5,8.2){$_{\times\ \Lambda'_k(q_\infty)}$}
\put(9.1,6.8){$_{D'_k}$}
\put(12.9,7.75){{\small -} $_0$}

\put(3,5.6){\vector(0,-1){0.6}}
\put(3.2,5.3){$\Lambda_{k+1}\circ\Lambda_k^{-1}$}
\put(11,5.6){\vector(0,-1){0.6}}
\put(11.2,5.3){$\Lambda'_k\circ{\Lambda'_k}^{-1}$}

\qbezier[10](3,4.5)(3,4.2)(3,3.9)
\qbezier[10](11,4.5)(11,4.2)(11,3.9)

\qbezier(1.8,3.2)(8,1.6)(1.8,0)
\put(3.5,1.45){{\small $\times$}}
\put(3.3,1.3){$_{(-1,\pri 0)}$}
\put(4.85,1.5){{\small -} $_0$}
\put(2,0.4){$_\H$}

\put(6.6,1.6){\vector(1,0){0.8}}
\put(6.8,1.8){$G$}

\qbezier(9.8,3.2)(16,1.6)(9.8,0)
\put(11.5,1.45){{\small $\times$}}
\put(11.3,1.3){$_{(-1,\pri 0)}$}
\put(12.85,1.5){{\small -} $_0$}
\put(10,0.4){$_\H$}

\end{picture}

\subsection{Dilations}
Let $(p_k)$ be a sequence of points of $D$ converging to $p_\infty\in\partial D$, and set $q_k=F(p_k)$. According to the boundary distance preserving property, $(q_k)$ converges to $q_\infty=F(p_\infty)\in\partial D'$.

\subsubsection{Choice of local coordinates} 
In an adapted local coordinate system $\Phi:U\to\R^{2n+2}$ about $p_{\infty}$ with $\Phi(p_\infty)=0$, we identify $p_{\infty}$ with 0 and $U$ with $\R^{2n+2}$. Moreover, we may assume that $\rho\circ\Phi^{-1}$ is bounded for the $\mathcal{C}^1$ norm, and:
\begin{itemize}
\item[\textbullet] $J(0)=J'(0)=\mathcal{J}_{st}$;
\item[\textbullet] $D=\{p\in\R^{2n+2}/\ \rho(p)<0\}$ and $T_0(\partial D)=\{x_0=0\}$, where the defining function $\rho$ can be expressed by
$$\rho(z)=\Re z_0+\Re\left(z_0\sum_{j\ge 1}(\rho_{j}z_j+\rho_{\bar{j}}\bar{z}_j)\right)+P(\pri z,\pri \bar{z})+\rho_\epsilon(z)$$
with $P$ being a real homogeneous polynomial of degree 2 and $\rho_\epsilon(z)=o(||z||^2)$;
\item[\textbullet] $D'=\{p\in \R^{2n+2}/\ \rho'(p)<0\}$ et $T_0(\partial D')=\{x_0=0\}$, where the defining function $\rho$ can be expressed by
$$\rho'(z)=\Re z_0+\Re\left(z_0\sum_{j\ge 1}(\rho'_{j}z_j+\rho'_{\bar{j}}\bar{z}_j)\right)+Q(\pri z,\pri \bar{z})+\rho'_\epsilon(z)$$
with $Q$ being a real homogeneous polynomial of degree 2 and $\rho'_\epsilon(z)=o(||z||^2)$.
\end{itemize}

\subsubsection{Centering}
Recall that for every neighborhood $V$ of 0, one can find some constant $\delta>0$ such that for all $p\in V\cap\partial D$, the closed ball of radius $\delta$ centered at $p-\delta \overrightarrow{n}_p$ is in $D\cup\{p\}$ (where $\overrightarrow{n}_p$ denotes the outer normal to $D$ at $p$). Hence $\dot{p}_k\in\partial D$ and $\dot{q}_k\in\partial D'$ such that 
$$\dist (p_k,\partial D)=||p_k-\dot{p}_k||=d_k\ \ \text{and}\ \ \dist (q_k,\partial D')=||q_k-\dot{q}_k||=d'_k$$
are uniquely defined for some sufficiently large $k$. 
Then there exists a rigid motion $\phi_k:\R^{2n+2}\to\R^{2n+2}$ with $\phi_k(\dot{p}_k)=0$ and $\phi_k(p_k)=(-d_k,0,\hdots,0)$, verifying:
\begin{itemize}
\item[\textbullet] the tangent space to $\partial(\phi_k(D))$ at 0 is $\{\Re z_0=0\}$ and the complex tangent space to $\partial(\phi_k(D))$ at 0 (for the induced almost complex structure $(\phi_k)_*J$) is $\{0\}\times\C^{n}$; 
\item[\textbullet] $\phi_k$ converges to the identity mapping on any compact subset of $\R^{2n+2}$ with respect to the $\mathcal{C}^2$ topology.
\end{itemize}

Consequently, $\dot{J}^k=(\phi_k)_*J$ converges to $J$ in the $\mathcal{C}^1$ sense on any compact subset, and is expressed by
\begin{equation}\label{formeJ}
\dot{J}^k(0)=\left(\begin{array}{cc}\dot{J}^k_{(1,1)}(0)&0_{2,2n}\\ \dot{J}^k_{(2,1)}(0)&\dot{J}^k_{(2,2)}(0)\end{array}\right).
\end{equation}
The sequence $\rho_k=\rho\circ\phi_k^{-1}$ also converges to $\rho$ at second order with respect to the compact-open topology:
$$\rho_k(z)=\rho\circ\phi_k^{-1}(z)=\tau_k\left(\Re z_0+\Re\left(z_0\sum_{j\ge 1}(\rho_{j}^kz_j+\rho_{\bar{j}}^k\bar{z}_j)\right)+P^k(\pri z,\pri\bar{z})+\rho^k_\epsilon(z)\right)$$
where $P^k$ is some real homogeneous polynomial of degree 2 and $\rho^k_\epsilon(z)=o(||z||^2)$ uniformly in $k$. Let us also define the inhomogeneous dilation
$\delta_k:(z_1,\hdots,z_n)\mapsto\left(\frac{z_1}{d_k},\frac{z_2}{\sqrt{d_k}},\hdots,\frac{z_n}{\sqrt{d_k}}\right)$, and set 
$$\Lambda_k=\delta_k\circ\phi_k\circ\Phi,\quad D_k=\Lambda_k(D),\quad r_k=\frac{1}{d_k\tau_k}\rho\circ\Lambda_k^{-1}\quad\mathrm{and}\quad J^k=(\Lambda_k)_*J.$$

We construct in the same way a rigid motion $\phi'_k:\R^{2n+2}\to\R^{2n+2}$, $\dot{J}'^k=(\phi'_k)_*J$, ${\rho'}_k=\rho'\circ{\phi'_k}^{-1}$ and the dilation $\delta'_k$. We also define $\Lambda'_k$, $D'_k$, $r'_k$ and $J'^k$. Finally, we set
$$F_k=\Lambda'_k\circ F\circ\Lambda_k^{-1}:D_k\to D'_k.$$

\subsubsection{Convergence}
\noindent After dilation, one gets
$d_kr_k(z)=d_k\left(\Re z_0+P^k(\pri z,\pri\bar{z})\right)+O(d_k\sqrt{d_k})$,
 which gives:

\begin{lm}\label{cvdomaines}
The sequence $(r_k)$ converges at second order to $\tilde{r}$ with respect to the compact-open topology, and $D_k$ converges in the sense of local Hausdorff set convergence to $\tilde{D}=\{z\in\R^{2n+2}/\ \tilde{r}(z)<0\}$, where 
$$\tilde{r}(z)=\Re z_0+P(\pri z,\pri\bar{z}).$$
There is a similar statement for $(r'_k)$ and $D'_k$.
\end{lm}

\begin{lm}
The sequence of almost complex structures $(J^k)$, respectively $(J'^k)$, converges on any compact subset to a model structure $\tilde{J}$, respectively $\tilde{J}'$, in the $\mathcal{C}^1$ sense .
\end{lm}

\begin{pr}
We follow \cite{Leebis}.
Writing almost complex structures as matrices, we have 
$$J(z)=J(0)+\left(\begin{array}{cc}A(z)&B(z)\\C(z)&D(z)\end{array}\right)=\left(\begin{array}{cc}\mathcal{J}^{(1)}_{st}+A(z)&B(z)\\C(z)&\mathcal{J}^{(n-1)}_{st}+D(z)\end{array}\right)$$
and
$$\dot{J}^k(z)=\dot{J}^k(0)+\left(\begin{array}{cc}\dot{A}^k(z)&\dot{B}^k(z)\\\dot{C}^k(z)&\dot{D}^k(z)\end{array}\right)=\left(\begin{array}{cc}\dot{J}^k_{(1,1)}(0)+\dot{A}^k(z)&\dot{B}^k(z)\\\dot{J}^k_{(2,1)}(0)+\dot{C}^k(z)&\dot{J}^k_{(2,2)}(0)+\dot{D}^k(z)\end{array}\right)$$
where $\dot{A}^k\to A,\ \dot{B}^k\to B,\ \dot{C}^k\to C$ and $\dot{D}^k\to D$ with respect to the $\mathcal{C}^1$ topology on any compact subset. Let us define
\begin{eqnarray*}
J^k(z)&=&\left(\begin{array}{cc}\frac{1}{d_k}I_2&0\\0&\frac{1}{\sqrt{d_k}}I_{2n}\end{array}\right)\,\dot{J}^k(\delta_k^{-1}(z))\,\left(\begin{array}{cc}d_kI_2&0\\0&\sqrt{d_k}I_{2n}\end{array}\right)\\
 &=&\left(\begin{array}{cc}\dot{J}^k_{(1,1)}+\dot{A}^k(\delta_k^{-1}(z))&\frac{1}{\sqrt{d_k}}\dot{B}^k(\delta_k^{-1}(z))\\ \sqrt{d_k}\dot{J}^k_{(2,1)}+\sqrt{d_k}\dot{C}^k(\delta_k^{-1}(z))&\dot{J}^k_{(2,2)}+\dot{D}^k(\delta_k^{-1}(z))\end{array}\right)\,.
\end{eqnarray*}
Since $\delta_k^{-1}$ converges uniformly to 0 and $\dot{J}^k$ converges uniformly to $J$ on any compact subset, it follows that
$$\dot{J}^k_{(1,1)}+\dot{A}^k(\delta_k^{-1}(z))\to \mathcal{J}^{(1)}_{st}$$
$$\sqrt{d_k}\dot{J}^k_{(2,1)}+\sqrt{d_k}\dot{C}^k(\delta_k^{-1}(z))\to 0$$
$$\dot{J}^k_{(2,2)}+\dot{D}^k(\delta_k^{-1}(z))\to \mathcal{J}^{(n)}_{st}$$
on any compact subset with respect to the $\mathcal{C}^1$ topology. Thus $\dot{B}^k(z)$ and $B(z)$ may be expressed by
$$\dot{B}^k(z)=\sum_{j=1}^n(B^k_{2j-1}x_j+B^k_{2j}y_j)+B^k_\epsilon(z)\quad\mathrm{and}\quad B(z)=\sum_{j=1}^n(B_{2j-1}x_j+B_{2j}y_j)+B_\epsilon(z),$$
where $B^k_j$ is a sequence of constant matrices which converges to $B_j$ as $k\to +\infty$, $B^k_\epsilon\to B_\epsilon$ in the $\mathcal{C}^1$ sense on any compact subset and $B^k_\epsilon=o(||z||)$ uniformly in $k$. Therefore,
\begin{eqnarray*}
\frac{1}{\sqrt{d_k}}\dot{B}^k(\delta_k^{-1}(z))
 &=&\sqrt{d_k}(B^k_1x_1+B^k_2y_1)+\sum_{j=2}^n(B^k_{2j-1}x_j+B^k_{2j}y_j)+\frac{1}{\sqrt{d_k}}B^k_\epsilon(d_kz_1,\sqrt{d_k}\pri z)\\
 &\to&\sum_{j=2}^n(B_{2j-1}x_j+B_{2j}y_j)\ \text{as}\ k\to +\infty.
\end{eqnarray*}
We obtain that $J^k$ converges on any compact subset of $\R^{2n+2}$ with respect to the $\mathcal{C}^1$-topology to $\tilde{J}$ defined as
$$\tilde{J}(z)=\left(\begin{array}{cc}\mathcal{J}_{st}^{(1)}&\tilde{B}(\pri z)\\0&\mathcal{J}_{st}^{(n)}\end{array}\right)\quad \text{where}\quad \tilde{B}(\pri z)=\sum_{j=2}^n(B_{2j-1}x_j+B_{2j}y_j).$$
\end{pr}




\begin{lm}
$(\tilde{D},\tilde{J})$ and $(\tilde{D}',\tilde{J}')$ are model domains.
\end{lm}

\begin{pr}
We recall the proof given in \cite{GSmodele} for the sake of completeness. Define $\tilde{r}_k=\rho\circ\delta_k^{-1}$ and $\tilde{J}^k={\delta_k}_*J$. As for $\dot{r}_k$ and $\dot{J}^k$, one can show that $\tilde{r}_k/d_k$ converges to $\tilde{r}$ at second order with respect to the compact-open topology, and $\tilde{J}^k$ converges to $\tilde{J}$ in the $\mathcal{C}^1$ sense on any compact subset. Consequently, for any $v$, $\mathcal{L}^{\tilde{J}^k}_0\left(\frac{\tilde{r}_k}{d_k}\right)(v)\xrightarrow[k\to +\infty]{}\mathcal{L}^{\tilde{J}}_0\tilde{r}(v)$.
By the invariance of the Levi form, we get $\mathcal{L}^J_0\rho(v)=\mathcal{L}^{\tilde{J}^k}_0\tilde{r}_k(d\delta_k(v))$. Since $\tilde{J}^k(0)=\mathcal{J}_{st}$, any complex tangent vector to the domain defined by $\tilde{r}_k$ is of the form $(0,v')$, and $d\delta_k(v)=v/\sqrt{d_k}$. For such a $v$,
$$\mathcal{L}^J_0\rho(v)=\mathcal{L}^{\tilde{J}^k}_0\tilde{r}_k(d\delta_k(v))=\mathcal{L}^{\tilde{J}^k}_0\tilde{r}_k(v/\sqrt{d_k})=\mathcal{L}^{\tilde{J}^k}_0\left(\frac{\tilde{r}_k}{d_k}\right)(v).$$
Passing to the limit, we obtain that $\mathcal{L}^{\tilde{J}}_0\tilde{r}(v)>0$ for every $v$ in the complex tangent space to $\tilde{D}$ at 0.
\end{pr}




\begin{lm}
The sequence $(F_k)$ admits a subsequence that converges at first order with respect to the compact-open topology to a map $\tilde{F}$ defined on $\tilde{D}$ and valued in $\overline{\tilde{D}'}$. The map $\tilde{F}$ is $(\tilde{J},\tilde{J}')$-holomorphic and verifies $\tilde{F}(-1,\pri 0)=(-1,\pri 0)$.
\end{lm}

\begin{pr}
Suppose that $K$ is some compact subset in $\tilde{D}$. 
Let us notice that, in order to get the existence of the desired subsequence, we just need to prove that $(F_k)$ is bounded in the $\mathcal{C}^0$-norm on $K$. Indeed, cover $K$ by small bidiscs, and consider two transversal foliations by $J$-holomorphic curves on every bidisc. Such a foliation is a small perturbation of the foliation by complex lines since $J$ is a small perturbation of the standard structure (see \cite{NW}). The restriction of $F_k$ on every such curve is bounded in the $\mathcal{C}^0$-norm, hence it is bounded in the $\mathcal{C}^1$-norm by the elliptic estimates \cite{S}. Since the bounds are uniform with respect to curves, the sequence $(F_k)$ is bounded in the $\mathcal{C}^1$-norm on $K$.

Now we follow \cite{Leebis}. First, we recall that there is some $\alpha>0$ such that for any sufficiently large $k$ and $r\in [0;1[$, and any $J^k$-holomorphic disc $h:\D\to D_k\cap U$ such that $h(0)\in Q(0,\alpha)$, there is some positive constant $C_r$ such that $h(\D_r)\subset Q(0,C_r\alpha)$. 

For every $p\in\tilde{D}$, there is some neighborhood $U_p$ of $p$, and some family $\mathcal{H}_p$ of pseudo-holomorphic discs centered at $p$, such that  $U_p\subset\bigcup_{h\in\mathcal{H}_p}h(\Delta_{r(p)})$ (see \cite{ref7Lee, ref27Lee, ref37Lee}). Hence one can find a finite covering $\{U_{t_j}\}_{j=0,\hdots,m}$ of $K$ such that $t_0=(-1,\pri 0)$ and $U_{t_j}\cap U_{t_{j+1}}\not=\emptyset$. Set $r=\text{max}\{r(t_j)\}$.
Since ${\delta'}_k^{-1}\circ F_k(-1,\pri 0)=(-d'_k,\pri 0)\in Q(0,d'_k)$, we have ${\delta'}_k^{-1}\circ {F}_k\circ h(\Delta_r)\subset Q(0,C_rd'_k)$ for all $h\in\mathcal{H}_{t_0}$, and 
$${\delta'}_k^{-1}\circ {F}_k(U_{t_0})\subset Q(0,C_rd'_k).$$ 
For all $h\in\mathcal{H}_{t_1}$, there exists $\omega\in\Delta_r$ such that $h(\omega)\in U_{t_0}\cap U_{t_1}$. The pseudo-holomorphic disc $g:\zeta\mapsto h\left(\frac{\zeta+\omega}{1+\bar{\omega}\zeta}\right)$ verifies $g(0)\in Q(0,C_rd'_k)$ and $g(\omega)=h(0)$. Thus ${\delta'}_k^{-1}\circ {F}_k(t_1)\in Q(0,C_r^2d'_k)$ and
${\delta'}_k^{-1}\circ {F}_k(U_{t_1})\subset Q(0,C_r^2d'_k)$. 
Iterating this process, one may obtain ${\delta'}_k^{-1}\circ {F}_k(U_{t_m})\subset Q(0,C_r^{2m+1}d'_k)$. Whence ${F}_k(K)\subset\delta'_k(Q(0,C_Kd_k'))=Q(0,C_K)$. 

Let $({F}_{k'})$ be a subsequence converging at first order with respect to the compact-open topology to some map $\tilde{F}:\tilde{D}\to\overline{\tilde{D}'}$. Passing to the limit in the pseudo-holomorphy condition, we get that $\tilde{F}$ is $(\tilde{J},\tilde{J}')$-holomorphic.
\end{pr}

\noindent In the particular case $n=2$, one can also find a proof of this statement in \cite{CGS},  based on the method developed in \cite{Berteloot, BertelootC}.

\subsection{Properties of the limit map \boldmath $G$\unboldmath}
There exist simple model structures $\mathcal{J}$ and $\mathcal{J}'$ on $\H$, and pseudo-biholomorphisms $\Psi:\tilde{D}\to\H$ and $\Psi':\tilde{D}'\to\H$ fixing $(-1,\pri 0)$, continuous and one-to-one to the boundary. Let us define $G=\Psi'\circ\tilde{F}\circ\Psi^{-1}$. By construction, $G:\H\to\overline{\H}$ is $(\mathcal{J},\mathcal{J}')$-holomorphic and fixes $(-1,\pri 0)$.

If the almost complex structure $\mathcal{J}$ (resp. $\mathcal{J}')$ is integrable, we can prescribe $\mathcal{J}=\mathcal{J}_{st}$ (resp. $\mathcal{J}'=\mathcal{J}_{st}$). 

\subsubsection{Boundary distance preserving property}
\begin{lm}
For any bounded subset $K$ in $\H$, there exist some constants $C_K,C_K'>0$ such that for all $p\in K$, $$C_K\le\frac{\dist(G(p),\partial\H)}{\dist(p,\partial\H)}\le C_K'.$$
In particular, $G$ takes its values in $\H$ (and not only $\overline{\H}$), and admits a locally Hölder $1/2$-continuous extension to $\overline{\H}$ verifying $G(\partial\H)\subset\partial\H$.
\label{conservationdistlim}
\end{lm}

\begin{pr}
The proof of Proposition \ref{conservationdist} gives two constants $c,c'>0$ such that for any $p\in D$, 
$$|\rho'(F(p))|\ge c\,\dist(p,\partial D)\quad\mathrm{and}\quad |\rho(p)|\ge c'\,\dist(F(p),\partial D').$$
Since ${F}_k=\Lambda'_k\circ F\circ\Lambda_k^{-1}$, we have for all $p\in {D}_k=\Lambda_k(D)$,
\begin{equation}\label{estk}
\begin{array}{l}
c\,\dist(\Lambda_k^{-1}(p),\partial D)\le|\rho'\circ{\Lambda'}_{k}^{-1}({F}_k(p))|=d'_k\tau'_k|{r'}_k({F}_k(p))|\\
c'\,\dist(F\circ\Lambda_k^{-1}(p),\partial D')\le|\rho\circ\Lambda_{k}^{-1}(p)|=d_k\tau_k|{r}_k(p)|.
\end{array}
\end{equation}
Let $q$ be a point of the boundary $\partial D$. Then
$$||\Lambda_k^{-1}(p)-q||\ge\frac{1}{\underset{\overline{D}}{\text{Max}}||d\phi_k||}\times||\delta_k^{-1}(p-\Lambda_k(q))||\ge \frac{d_k}{\underset{\overline{D}}{\text{Max}}||d\phi_k||}\times||p-\Lambda_k(q)||.$$
Therefore, we obtain by (\ref{estk}):
$$c_k\,\dist(p,\partial\tilde{D}_k)\le|\tilde{r}'_k({F}_k(p))|\quad\mathrm{where}\quad c_k=c\,\frac{d_k}{d'_k\tau'_k\underset{\overline{D}}{\text{Max}}||d\phi_k||};$$
$$c'_k\,\dist(F_k(p),\partial\tilde{D}'_k)\le|\tilde{r}_k(p)|\quad\mathrm{where}\quad c'_k=c'\,\frac{d'_k}{d_k\tau_k\underset{\overline{D}}{\text{Max}}||d\phi'_k||}\,.$$
According to Proposition \ref{conservationdist}, $\frac{d'_k}{d_k}=\frac{\text{dist}(F(p_k),\partial D')}{\text{dist}(p_k,\partial D)}$ is bounded between to positive constants. Passing to the limit, we get some constants $C,C'>0$ such that for any $p\in\tilde{D}$ (and hence $p\in\tilde{D}_k$ for some sufficiently large $k$),
\begin{equation}\label{distFtilde}
\begin{array}{l}
c\,\dist(\tilde{F}(p),\partial\tilde{D}')\le |\tilde{r}(p)|\\
c'\,\dist(p,\partial\tilde{D})\le |\tilde{r}'(\tilde{F}(p))|.
\end{array}
\end{equation}
Applying the diffeomorphisms $\Psi'$ and $\Psi^{-1}$, we obtain the two desired inequalities.
The same arguments as in the proof of Corollary \ref{cor:prolC1/2} give the locally Hölder $1/2$-continuous extension to $\overline{\H}$.
\end{pr}

\begin{cor}
Writing $G=(G_0,\pri G)$, one has: $\Re(G_0(t,\pri 0))\xrightarrow[t\in\R,\ t\to -\infty]{}-\infty$.
\label{valeurinfini}
\end{cor}

\begin{pr}
Since $r(z)=\Re z_0+||\pri z||^2\ge\Re z_0$, we just have to show that $r(G(t,\pri 0))\to -\infty$. Moreover, for any $t\in\R^-$,
$$r(G(t,\pri 0))=\tilde{r}'(\tilde{F}(\Psi^{-1}(t,\pri 0)))=\tilde{r}'(\tilde{F}(t,\pri 0))\ge c\,\dist((t,\pri 0),\partial\tilde{D})$$
according to (\ref{distFtilde}). So it suffices to get $\dist((t,\pri 0),\partial\tilde{D})\xrightarrow[t\in\R,\ t\to -\infty]{}\infty$. The domain $\tilde{D}$ is defined by $0=\Re z_0+P(\pri z,\pri\bar{z})$, where $P$ is a real homogeneous polynomial of degree two. Let $\gamma>0$ be such that for any $\pri z\in\C^n,\ |P(\pri z,\pri\bar{z})|\le\gamma||\pri z||^2$. 
One may easily verify that for any $z=(z_0,\pri z)\in\C^{n+1}$ such that $||(t,\pri 0)-z||<\sqrt{\frac{|t|}{1+\gamma}}$, we have $\tilde{r}(z)<0$
as soon as $\frac{|t|}{1+\gamma}\ge 1$. Whence, for $t$ sufficiently large, $\dist((t,\pri 0),\partial\tilde{D})\ge \sqrt{\frac{|t|}{1+\gamma}}$.
\end{pr}

\subsubsection{Studying the Jacobian}
In order to simplify the notations, we assume from now  that $\tilde{D}=\tilde{D}'=\H$ and $\Psi=\Psi'=id$.

\begin{lm}
There exist some constants $0<\alpha\le\beta<\infty$ such that for all $p\in\H$,
$$\alpha|\mathrm{Jac}_pG|\le\mathrm{liminf}\,\left|\mathrm{Jac}_{\Lambda_k^{-1}\circ\delta_k^{-1}(p)}F\right|\le\mathrm{limsup}\,\left|\mathrm{Jac}_{\Lambda_k^{-1}\circ\delta_k^{-1}(p)}F\right|\le\beta|\mathrm{Jac}_pG|.$$
\label{jacobienlim}
\end{lm}

\begin{pr}
Pick $p\in\H$. Then $p\in D_k$ for some sufficiently large $k$, and
$$d(F_k)_p=d(\delta'_k)\circ d(\Lambda'_k)\circ dF_{\Lambda_k^{-1}\circ\delta_k^{-1}(p)}\circ d\Lambda_k^{-1}\circ d\delta_k^{-1}.$$
The rigid motion $\Lambda_k$ converges to the identity mapping. Hence, taking the determinant in the previous equality, we get  
\begin{eqnarray}\label{egalitejac}
\mathrm{Jac}_{\Lambda_k^{-1}\circ\delta_k^{-1}(p)}F=\mu_k\mathrm{Jac}_pF_k,
\end{eqnarray}
where $\mu_k$ only depends on $k$. Moreover, $\mu_k\underset{k\to +\infty}{\sim}(d'_k/d_k)^{n+1}$ remains bounded between two positive constants by the boundary distance preserving propertie of $F$. Since $\mathrm{Jac}_pF_k=\mathrm{det}(d(F_k)_p)\xrightarrow[k\to +\infty]{}\mathrm{Jac}_pG$, we conclude the proof by passing to the limit in (\ref{egalitejac}).
\end{pr}

\begin{lm}\label{chgtP}
For any sequence $\mathcal{P}=(p_k)_k$ of points of $D$ converging to $p_\infty\in\partial D$, the sequence $(\mathrm{Jac}_{p_k}F)_k$ is bounded. If $\mathrm{Jac}_{p_k}F\underset{+\infty}{\to} 0$, then for any other sequence $(p'_k)_k$ of points of $D$ converging to $p_\infty$, we have $\mathrm{Jac}_{p'_k}F\underset{+\infty}{\to} 0$.
\end{lm}

\begin{pr}
Denote by $G^{\mathcal{P}}$ the limit map obtained by applying the scaling method to the sequence $\mathcal{P}$. Lemma \ref{jacobienlim} taken at $p=(-1,\pri 0)\in\H$ gives:
$$\alpha|\mathrm{Jac}_{(-1,\pri 0)}G^\mathcal{P}|\le\mathrm{liminf}\,|\mathrm{Jac}_{p_k}F|\le\mathrm{limsup}\,|\mathrm{Jac}_{p_k}F|\le\beta|\mathrm{Jac}_{(-1,\pri 0)}G^\mathcal{P}|.$$
This implies that the sequence $(\mathrm{Jac}_{p_k}F)_k$ is bounded.

Assume that $\mathrm{Jac}_{p_k}F\to 0$, and let $(p'_k)_k$ be another sequence converging to $p_\infty$. Set $p''_{2k}=p_k$ and $p''_{2k+1}=p'_k$. For any cluster point $\lambda$ of $(\mathrm{Jac}_{p'_k}F)$, the sequence $(\mathrm{Jac}_{p''_k}F)$ has at least 0 and $\lambda$ as cluster points. The scaling method applyed to the sequence $\mathcal{P}''=(p''_k)$ and Lemma \ref{jacobienlim} taken at $p=(-1,\pri 0)\in\H$ show: 
$$\alpha''|\mathrm{Jac}_{(-1,\pri 0)}G^{\mathcal{P}''}|\le 0\le|\lambda|\le\beta''|\mathrm{Jac}_{(-1,\pri 0)}G^{\mathcal{P}''}|.$$
Hence $\mathrm{Jac}_{(-1,\pri 0)}G^\mathcal{P''}=0$ and $\lambda=0$.
\end{pr}

\begin{lm}\label{jacobienG'}
Let $\mathcal{P}=(p_k)$ be a sequence of points of $D$ converging to $p_\infty\in\partial D$. Then the Jacobian of $G=G^{\mathcal{P}}$ does not vanish in $\H$.
\end{lm}

\begin{pr}
We may assume $p_\infty=0$.\\
\textbullet\ First, we show that if the Jacobian of $G$ vanishes at some point $p\in\H$, then it vanishes identically in $\H$.
The scaling method applyed to $\mathcal{P}$ and Lemma \ref{jacobienlim} taken at $p$ give  a sequence $(p''_k)$, where $p''_k=\Lambda_k^{-1}\circ\delta_k^{-1}(p)$, such that $\mathrm{Jac}_{p''_k}F\to 0$. Then, for any $p'\in\H$, we get
\begin{equation}\label{inegalitejac}
\alpha|\mathrm{Jac}_{p'}G|\le\mathrm{liminf}\,\left|\mathrm{Jac}_{\Lambda_k^{-1}\circ\delta_k^{-1}(p')}F\right|\le\mathrm{limsup}\,\left|\mathrm{Jac}_{\Lambda_k^{-1}\circ\delta_k^{-1}(p')}F\right|\le\beta|\mathrm{Jac}_{p'}G|.
\end{equation}
Hence $\mathrm{Jac}_{p'}G=0$ if and only if $\mathrm{Jac}_{\Lambda_k^{-1}\circ\delta_k^{-1}(p')}F\to 0$. According to Lemma \ref{chgtP}, it only remains to prove that the sequence $(p'_k)$ converges to 0, with $p'_k=\Lambda_k^{-1}\circ\delta_k^{-1}(p')$.
But 
$$\Lambda_k^{-1}\circ\delta_k^{-1}(p')\xrightarrow[k\to +\infty]{}0,$$
which gives the statement.\\
\textbullet\ Suppose by contradiction that the Jacobian of $G$ is identically zero in $\H$. There exist a neighborhood $U$ of 0, a constant $\delta>0$ and a function $\vp$ continuous on $\overline{U\cap\H}$, strictly $J$-plurisubharmonic on $U\cap\H$, such that:
\begin{equation}\label{ineq}
\forall z\in U\cap\H,\ \vp(z)< -\delta ||z||^2.
\end{equation}
Let us fix $\ve>0$ such that $\overline{\B(0,\sqrt{\ve/\delta})}\subset U$, and set $H^{\ve}=\{z\in U\cap\H/\ \vp(z)>-\ve\}$. Then $\overline{H^{\ve}}\subset U$ by (\ref{ineq}). 
By hypothesis, the maximal rank $r_0$ of $dG$ on $U\cap \H$ is less than $2n+1$. Moreover, according to Lemma \ref{conservationdistlim}, $G(U\cap\H)\subset U\cap\H$ and its continuous extension verifies $G(\partial\H)\subset\partial\H$. Hence $G$ is non-constant and $r_0>0$.

As in the proof of Lemma \ref{lm:valeursregulieres}, we obtain the existence of some $q\in G(U\cap \H)$ such that $N=G^{-1}(\{q\})$ is an almost complex submanifold of (real) dimension $2n+2-r_0$ in $U\cap\H$. The continuous function $\vp$ reaches its maximum on the compact set $N\cap\overline{H^{\ve}}$ at some point $p_0$. If $p_0\in N\cap H^{\ve}$, there exists a pseudo-holomorphic disc $h$ in the open subset $N\cap H^{\ve}$ of $N$, centered at $p_0$ and such that $\frac{\partial h}{\partial x}(0)\not=0$ \cite{NW}, but this is impossible by the maximum principle. Thus $p_0\in N\cap\partial H^{\ve}$. Because of the continuity of $\vp$, 
$$\partial H^{\ve}=(U\cap\partial\H)\cup\{z\in U\cap\H/\ \vp(z)=-\ve\}.$$ 
The boundary distance preserving property of $G$ implies that $N$ does not intersect $\partial\H$. Consequently, $\vp(p_0)=-\ve$, and $\underset{N\cap\overline{H^{\ve}}}{\mathrm{Max}}\vp=-\ve$. Hence $\vp$ is constant equal to $-\ve$ on $N\cap\overline{H^{\ve}}$, which contradicts the strict plurisubharmonicity of $\vp$.
\end{pr}

\subsubsection{Computation of $\frac{\partial G_0}{\partial z_0}$}
\begin{lm}\label{diffen10}
For all $z\in\H$, $\displaystyle\frac{\partial G_0}{\partial z_0}(z)=1$.
\end{lm}

\begin{pr}
First, notice that by Lemma \ref{jacobienG'}, the almost complex structures $\mathcal{J}$ and $\mathcal{J}'$ are both integrable either non-integrable.\\
\textbullet\ Suppose that $\mathcal{J}$ and $\mathcal{J}'$ are both integrable. Then, we have seen that we can assume $\mathcal{J}=\mathcal{J}'=\mathcal{J}_{st}$. The map $G:\H\to\H$ is thus holomorphic for the standard structure, and admits a continuous extension to the boundary (Lemma \ref{conservationdistlim}). Denote by $\Phi$ the biholomorphism (for the standard structure) from $\H$ to the unit ball $\B$ of $\C^{n+1}$ defined by 
$$\Phi(z_0,\pri z)\mapsto\left(\frac{z_0+1}{z_0-1},\frac{1}{1-z_0}\pri z\right).$$
It extends to a homeomorphism from $\overline{\H}$ to $\overline{\B}$ by defining $\Phi(\infty)=(1,\pri 0)$ and $\Phi^{-1}(1,\pri 0)=\infty$.

The map $\tilde{G}=\Phi\circ G\circ \Phi^{-1}$ from $\B$ to $\B$ is holomorphic, continuous up to $S^*=\partial\B\setminus\{(1,\pri 0)\}$. Moreover, $\tilde{G}(S^*)\subset\partial\B$. Such a map is an automorphism of the ball (see \cite{PinchukT}, Proposition 2.3).
We have $\tilde{G}(0)=0$, and for all $u\in[0;1[$:
$$\tilde{G}(u,\pri 0)=\Phi\circ G\left(\frac{u+1}{u-1},\pri 0\right)=\left(\frac{Z_0+1}{Z_0-1},\frac{\sqrt{1}}{1-Z_0}\pri Z\right)$$
where $G(\frac{u+1}{u-1},\pri 0)=(Z_0,\pri Z)$. Corollary \ref{valeurinfini} implies that as $u$ tends to $1^-$, the real part of $Z_0$ tends to $-\infty$ and so
$\Re\frac{Z_0+1}{Z_0-1}\to 1$.
Since the image of $\Phi$ is the unit ball, necessarily $\tilde{G}(u,\pri 0)\xrightarrow[u\in[0;1[,\,u\to 1]{}(1,\pri 0)$. According to \cite{Co} (p. 467), we get $\tilde{G}_0\equiv id$. Hence $G_0(z_0,\pri z)=z_0$ for every $z\in\H$.

\noindent\textbullet\ Suppose that $\mathcal{J}$ et $\mathcal{J}'$ are both non-integrable.
In this case, we get by (\ref{FormeCasNonInt}) that 
$$G(z_0,\pri z)=(cz_0+f_1(\pri z)+if_2(\pri z),\pri G(\pri z)),$$
where $c\not=0$ is a real constant, and $f_1$, $f_2$ are real valued. Since the map $G$ is continuous to the boudary and verifies $G(\partial\H)\subset\partial\H$, we have for any $\pri z\in\C^n$:
$$\Re z_0+||\pri z||^2=0\Longrightarrow c\,\Re(z_0)+f_1(\pri z)+||\pri G(\pri z)||^2=0.$$
As a consequence, $f_1(\pri z)=c||\pri z||^2-||\pri G(\pri z)||^2$, $f_1(\pri 0)=||\pri G(\pri 0)||^2$ and 
$$(-1,\pri 0)=G(-1,\pri 0)=(-c+f_1(\pri 0)+if_2(\pri 0),\pri G(\pri 0))\Longrightarrow \left\{\begin{array}{l}\pri G(\pri 0)=\pri 0\\ c=1\end{array}\right..$$
Thus $\displaystyle\frac{\partial G}{\partial z_0}(z)=1$ for all $z\in\H$.
\end{pr}

\section{Proofs of the theorems}
\subsection{Behavior near the boundary}
The different properties of $G$ proved in the previous section give us some informations on $F$ near the boundary. In particular, Lemma \ref{jacobienG'} implies that there is no sequence $(p_k)$ of points of $D$ converging to some point of the boundary such that $\mathrm{Jac}_{p_k}F\to 0$. This proves Theorem \ref{theo:jacobiencompact}.\\ 

As a consequence, the map $F$ is locally pseudo-biholomorphic out of some compact set. Whence, in order to study the regularity near the boundary, it suffices to understand the pseudo-biholomorphic case.
We also know by \cite{CGS}, Proposition 3.5, some precise estimates of the Kobayashi metric, which give the asymtotic behavior of the differential according to the directions. We begin with fixing notations. 
Consider the vector fields
$$v^0=(\partial\rho/\partial x_0)\partial/\partial y_0-(\partial\rho/\partial y_0)\partial/\partial x_0$$
and
$$v^j=(\partial\rho/\partial x_0)\partial/\partial x_j-(\partial\rho/\partial x_j)\partial/\partial x_0\quad\mathrm{for}\quad j=1,\hdots,n.$$
Restricting if necessary the neighborhood $U$ of 0 on which we are working, the vector fields defined by $X^j=v^j-iJv^j$, $1\le j\le n$ form a basis of the $J$-complex tangent space to $\{\rho=\rho(z)\}$ at any $z\in U$. Moreover, if $X^0=v^0-iJv^0$, then the family $X=(X^0,X^1,\hdots,X^n)$ forms a basis of $(1,0)$ vector fields on $U$. Similarly, we construct a basis $X'=(X'^0,X'^1,\hdots,X'^n)$ of $(1,0)$ vector fields on $U'$ such that $(X'^1(w),\hdots,X'^n(w))$ defines a basis of the $J'$-complex tangent space to $\{\rho'=\rho'(w)\}$ at any $w\in U'$. In the sequel, we will denote by $A(p_k)$ the matrix of the map $dF_{p_k}$ with respect to $X(p_k)$ and $X'(F(p_k))$.

\begin{prop}\label{estimationdiff}
{\bf \cite{CGS}} The matrix $A(p_k)$ satisfies the following estimates: 
$$A(p_k)=\left(\begin{array}{cc}
O_{1,1}(1)& O_{1,n}(\dist(p_k,\partial D)^{1/2})\\
O_{n,1}(\dist(p_k,\partial D)^{-1/2})& O_{n,n}(1)
\end{array}\right)\,.$$
\end{prop}
\begin{rem}
The asymptotic behaviour of $A(p_k)$ depends only on the distance from $p_k$ to $\partial D$, not on the choice of the sequence $(p_k)_k$.
\end{rem}

In the case of a biholomorphism, one gets immediately a similar estimate for $(dF_{p_k})^{-1}=d(F^{-1})_{F(p_k)}$. For a proper map, the control of the inverse matrix comes from Proposition \ref{estimationdiff} and the control of the Jacobian.

\begin{prop}
The matrix $A(p_k)$ is invertible, and its inverse verifies the following estimates:
$$A(p_k)^{-1}=\left(\begin{array}{cc}
O_{1,1}(1)& O_{1,n}(\dist(p_k,\partial D)^{1/2})\\
O_{n,1}(\dist(p_k,\partial D)^{-1/2})& O_{n,n}(1)
\end{array}\right)\,.$$
\label{estimationinversediff}
\end{prop}

\begin{dem}
The formula $A^{-1}=\frac{1}{\mathrm{Jac}\,F}\times\tsp\mathrm{com}\,A$, together with Lemma \ref{chgtP} and Theorem \ref{theo:jacobiencompact}, shows that we just need to get estimates for the matrix $B=\tsp\mathrm{com}\,A$. The determinant extracted from $A$ and appearing in the entry $B_{i,j}$ can be expressed by developping along the lign number 0 and/or the column number 0 of $A$. This gives the statement according to Proposition \ref{estimationdiff}.
\end{dem}

Lemma \ref{diffen10} provides, exactly as in Proposition 4.5 in \cite{GSmodele}, more informations on the entry $(0,0)$ of the matrix $A(p_k)$. Note first that $X$ and $X'$, and hence the matrix $A(p_k)$, were normalized by the condition $J(p_{\infty})=\mathcal{J}_{st}$. So let us note $A(p_\infty,p_k)$ in place of $A(p_k)$.

\begin{prop}
The entry $(0,0)$ of the matrix $A$ verifies the following properties:
\begin{itemize}
\item every cluster point of the function $z\mapsto A_{0,0}(p,z)$ is real when $z$ tends to $p\in\partial D$;
\item given $z\in D$, let $p\in\partial D$ be such that $\dist(z,\partial D)=||p-z||$; then, there exists a constant $A>0$, independent of $z\in D$, such that $|A_{(0,0)}(p,z)|\ge A$.
\end{itemize}
\label{limreellenonnulle}
\end{prop}

By means of Theorem \ref{theo:jacobiencompact} and Propositions \ref{estimationinversediff} and \ref{limreellenonnulle}, we may use the arguments of the proof of Theorem 0.1 in \cite{GSmodele}. As a consequence, we get:

\begin{theo}
Let $(D,J)$ and $(D',J')$ be strictly pseudoconvex regions of the same dimension. Then, every proper pseudoholomorphic map from $D$ to $D'$ has a $\mathcal{C}^1$-extension to $\overline{D}$.
\label{regularite}
\end{theo}

\subsection{Proof of Theorem \ref{theo:estimationpropre} and Corollary \ref{cor:CNS}}
We suppose that the conditions of Theorem \ref{theo:estimationpropre} are satisfied. Since $F$ is a local biholomorphism near the boundary, we can apply Corollary 1 of \cite{3} to the map $(F,\tsp(dF)^{-1})$ with $N=N^*M$ et $N'=N^*M'$. We obtain that the map $(F,\tsp(dF)^{-1})$ is locally of class $\mathcal{C}^{t_1-1}$, where 
$t_1=\mathrm{min}\,(r-1,r')$.
We similarly get that the map $(F^{-1},\tsp(dF))$ (which is well-defined near the boundary) is locally of class $\mathcal{C}^{t_2-1}$, where 
$t_2=\mathrm{min}\,(r'-1,r)$. 

Thus $F$ is of class $\mathcal{C}^{s}$, where $s=\mathrm{max}\,(t_1,t_2)$. This gives $s=\mathrm{max}\,(r-1,r'-1)$ if $|r'-r|<1$, and $s=\mathrm{min}\,(r,r')$ if $|r'-r|\ge 1$. 
Moreover, 
$$||F||_{\mathcal{C}^{s-1}(\bar{D})}\le c(s)||(F,\tsp(dF)^{-1})||_\infty\left(1+\frac{c'}{\sqrt{\lambda_{N'}^{J'}}}\right).$$
This concludes the proof of Theorem \ref{theo:estimationpropre}.

Corollary \ref{cor:CNS} follows immediately, since the almost complex structure $J'$ is defined near the boundary by $J'_q=dF_q\circ J\circ (dF_q)^{-1}$.


\address
\email


\end{document}